\documentclass[nosumlimits,twoside,12pt]{amsart}
\usepackage{amssymb}
\usepackage{amsmath}
\usepackage[bookmarksnumbered,plainpages,hyperindex]{hyperref}
\usepackage{amsfonts}
\usepackage[usenames]{color}
\usepackage{float}
\usepackage{graphicx}
\newtheorem{theorem}{\sc Theorem}[section]
\newtheorem{proposition}[theorem]{\sc Proposition}
\newtheorem{lemma}[theorem]{\sc Lemma}
\newtheorem{corollary}[theorem]{\sc Corollary}
\theoremstyle{definition}
\newtheorem{definition}[theorem]{\sc Definition}

\theoremstyle{remark}
\newtheorem{remark}[theorem]{\sc Remark}

\def\ker{\mathrm{Ker\,}}

\def\ot{\otimes}

\def\H{\mathbb{H}}

\def\F{\mathbf{F}}

\input{xy}
\xyoption{all}
\newcommand{\diagBradedl}{\xymatrix@R=30pt@C=3pt{
&B\otimes B\ot B\otimes B \ar[rr]^{ B\ot c\ot B} &&B\otimes B\ot
B\otimes B
\ar[dr]|{ m\ot m}&\\
B\ot B\ar[ur]|{\Delta\ot \Delta } \ar[drr]|{m} &&&&B\otimes B
\\&&B\ar[urr]|{\Delta }&&&}}
\newcommand{\diagBradedr}{\xymatrix@C=40pt{
B\otimes B \ar[d]_{m} \ar[r]^{\varepsilon \otimes \varepsilon} & \mathbf{1}\ot \mathbf{1} \ar[d]^{m_{\mathbf{1}}} \\
B \ar[r]^{\varepsilon} & \mathbf{1}   }}
\newcommand{\diagCosepl}{\xymatrix@R=15pt@C=50pt{
{^C\mathcal{M}^C} \ar[d]_{\mathbb{T}} \ar[r]^{F'} & {^{C'}\mathcal{M'}^{C'}} \ar[d]^{\mathbb{T}'} \\
\mathcal{M}  \ar[r]_{F} & \mathcal{M'}   }}
\newcommand{\diagCosepr}{\xymatrix@R=15pt@C=50pt{
{^C\mathcal{M}^C} \ar[r]^{F'} & {^{C'}\mathcal{M'}^{C'}}  \\
\mathcal{M}  \ar[u]^{\mathbb{H}} \ar[r]_{F} & \mathcal{M'}
\ar[u]_{\mathbb{H}'}  }}
\newcommand{\diagFSl}{\xymatrix@R=15pt@C=50pt{
{^B\mathfrak{M}_B^B} \ar[d]_{\mathbb{T}} \ar[r]^{F'} & {^{B}\mathfrak{M}^{B}} \ar[d]^{\mathbb{T}'} \\
\mathfrak{M}_B  \ar[r]_{F} & \mathfrak{M}   }}
\newcommand{\diagFSr}{\xymatrix@R=15pt@C=50pt{
{^B\mathfrak{M}^B_B} \ar[r]^{F'} & {^{B}\mathfrak{M}^{B}}  \\
\mathfrak{M}_B  \ar[u]^{\mathbb{H}} \ar[r]_{F} & \mathfrak{M}
\ar[u]_{\mathbb{H}'}  }}

\def\id{\rm{id}}
\newcommand{\diagQcatA}{\xymatrix{ M\ar[r]^f \ar[d] & N\ar[d] \\M' \ar[r]^g & N'}}
\newcommand{\diagQcatB}{\xymatrix{M\ar[r]^f & N}}
\newcommand{\diagQcatC}{\xymatrix{ M\ar[r]^f \ar[d]_{{\id}_M} & N\ar[d]^{{\id}_N} \\
M \ar[r]^f & N}}
\def\eclass#1{{\mathcal{E}_{#1}}}

 \newcounter{zlist}
  \newenvironment{zlist}{\begin{list}{(\arabic{zlist})}{
  \usecounter{zlist}\leftmargin2.5em\labelwidth2em\labelsep0.5em
  \topsep0.6ex
  \parsep0.3ex plus0.2ex minus0.1ex}}{\end{list}}

  \newcounter{blist}
  \newenvironment{blist}{\begin{list}{(\alph{blist})}{
  \usecounter{blist}\leftmargin2.5em\labelwidth2em\labelsep0.5em
  \topsep0.6ex 
  \parsep0.3ex plus0.2ex minus0.1ex}}{\end{list}}

  \newcounter{rlist}
  \newenvironment{rlist}{\begin{list}{(\roman{rlist})}{
  \usecounter{rlist}\leftmargin2.5em\labelwidth2em\labelsep0.5em
  \topsep0.6ex 
  \parsep0.3ex plus0.2ex minus0.1ex}}{\end{list}}

\begin{document}
\title{Formally smooth bimodules}
\author{A.\ Ardizzoni}
\address{University of Ferrara, Department of Mathematics, Via Machiavelli 35,\newline \indent Ferrara,
I-44100, Italy} \email{rdzlsn@unife.it}
\urladdr{http://www.unife.it/utenti/alessandro.ardizzoni}
\author{Tomasz Brzezi\'nski}
\address{Department of Mathematics, University of Wales Swansea, Singleton
Park,\newline \indent Swansea SA2 8PP, U.K.}
\email{T.Brzezinski@swansea.ac.uk}
\urladdr{http://www-maths.swan.ac.uk/staff/tb}
\author{C.\ Menini}
\address{University of Ferrara, Department of Mathematics, Via Machiavelli 35,\newline \indent Ferrara,
I-44100, Italy} \email{men@unife.it}
\urladdr{http://www.unife.it/utenti/claudia.menini}
\thanks{This paper was written when A.\ Ardizzoni and C.\ Menini were members of
G.N.S.A.G.A.\ with partial financial support from M.I.U.R. The stay of
T.\ Brzezi\'nski  at University of Ferrara was
supported by I.N.D.A.M.}

\begin{abstract}
The notion of a formally smooth bimodule is introduced and its basic
properties are analyzed. In particular it is proven that a $B$-$A$ bimodule $M$
which is a generator left $B$-module is formally smooth if and only if the
$M$-Hochschild dimension of $B$ is at most one. It is also shown that modules
$M$ which are generators in the category $\sigma[M]$ of $M$-subgenerated
modules provide natural examples of formally smooth bimodules.

\end{abstract}
\maketitle

\section{Introduction}
The notion of a {\rm (formally) smooth algebra} was introduced in \cite{Sch:smo}. It has
been recognized in \cite{CunQui:alg} that smooth (or quasi-free) algebras can be
interpreted as functions on non-commutative nonsingular (smooth) affine
varieties or as
analogues of manifolds in non-commutative geometry. This point of view was then
developed further in \cite{KonRos:smo}, where an approach to {\em smooth
non-commutative geometry} was outlined. In
\cite{KonRos:non} this has given rise to
the introduction of  formally smooth objects, morphisms and functors as main building blocks of non-commutative algebraic
geometry. Following on, the non-commutative
geometric aspects of smooth algebras (or, more generally, {\em $R$-rings} or {\em
smooth algebra extensions}) such as tangent and cotangent bundles or
symplectic structures have been  discussed in \cite{CraEti:non}
(cf.\ \cite{vdB:dou}),
in the framework of {\em double derivations}. A general algebraic approach to formal
smoothness in monoidal abelian categories, including the cohomological aspects, has
been recently proposed in \cite{Ar1} and \cite{AMS1}.

The aim of this paper is to find a common ground for the notions of formal
smoothness which have attracted so much attention in recent literature. The basic
idea for this goes back to \cite{Su}, where it is observed that properties of
an extension of algebras, such as separability, can be encoded more generally
as properties of bimodules rather than algebra maps. We thus propose the
definition of a {\em formally smooth bimodule}, and show that this notion encodes
 smooth algebras and smooth extensions (which can be understood as
smooth algebras in monoidal category of bimodules). Furthermore we
show that a smooth bimodule can be interpreted as a smooth object
in the sense of \cite{KonRos:non}.  The definition of a smooth
bimodule is presented within the framework of relative homological
algebra, making specific use of tools recently developed in
\cite{Ar1}, and, in particular, developing the {\em
bimodule-relative cohomology}. With these tools we show that
separable bimodules can be understood as (non-commutative,
relative) ``bundles of points" (objects with zero
relative-Hochschild dimension), while the formally smooth
(generator) bimodules can be viewed as (non-commutative, relative)
``bundles of curves" or ``line bundles" (objects with
relative-Hochschild dimension at most one). On more
module-theoretic side, we show that given a left $B$-module $M$
with endomorphism ring $S$, $M$ is a separable $B$-$S$ bimodule if
and only if it is a generator of all left $B$-modules. On the
other hand, if ${M}$ is a generator in the category
$\sigma[M]$ of $M$-subgenerated left $B$-modules, then $M$ is a
formally smooth $B$-$S$ bimodule.

\subsection*{Module-theoretic conventions} By a ring we mean a unital associative
ring. ${}_B\mathcal{M}$, $\mathcal{M}_A$, ${}_B\mathcal{M}_A$ denote
categories of (unital) left $B$-modules, right $A$-modules and $B$-$A$
 bimodules. Morphisms in these categories are respectively denoted by
 ${}_B\mathrm{Hom}(-,-)$,  $\mathrm{Hom}_A(-,-)$ and  ${}_B\mathrm{Hom}_A(-,-)$.
 For a $B$-$A$  bimodule $M$ we often write ${}_BM$, $M_A$, ${}_BM_A$ to
 indicate the ring  and module structures used. The arguments of left $B$-module
 maps are always written on the left. This induces a composition convention for the endomorphism ring $S:={}_B\mathrm{End}(M)$ of ${}_BM$, which makes $M$ a $B$-$S$ bimodule. Given
 bimodules ${}_BM_A$ and ${}_BN_T$, we view the abelian group
 ${}_B\mathrm{Hom}(M,N)$ as an $A$-$T$ bimodule with multiplications defined
 by
 $$
 (m)(aft) := (ma)\!ft, \qquad \mbox{for all } f\in {}_B\mathrm{Hom}(M,N),\; a\in A,\;
 m\in M,
 \; t\in T.
 $$
 For a $B$-$A$ bimodule $M$, ${}^*M$ denotes the dual $A$-$B$ bimodule
 ${}_B\mathrm{Hom}(M,B)$.

\section{Relative projectivity and separable functors}

\subsection{Relative projectivity and injectivity}
A convenient description and conceptual interpretation of formally smooth or
separable bimodules is provided by relative cohomology. In this introductory section
we recall the basic properties of relative derived functors.  Most of the material
presented here can be found in \cite[Chapter IX]{HS}.

Let $\mathfrak{C}$ be a
category and let $\mathcal{H}$ be a class of morphisms in $\mathfrak{C}$.
An object $P\in\mathfrak{C}$ is called $f$-\emph{projective}, where
$f:C_{1}\rightarrow C_{2}$ is a morphism in $\mathfrak{C}$, if
\[
\mathfrak{C}(P,f): {\mathfrak{C}}(P,C_{1}
)\rightarrow{\mathfrak{C}}(P,C_{2}), \qquad g\mapsto f\circ g
\]
is surjective. $P$ is said to be $\mathcal{H}$-\emph{projective} if it is
$f$-projective for every $f\in\mathcal{H}$.

The \emph{closure} $\overline{\mathcal{H}}$
of the class of morphisms $\mathcal{H}$ is defined by
\[
\overline{\mathcal{H}}:=\left\{  f\in\mathfrak{C}\mid \text{if an object } P\in
\mathfrak{C}\text{ is }
\mathcal{H}\text{-projective, then } P\text{ is }f\text{-projective}\right\}.
\]
Obviously, $\overline{\mathcal{H}}$ contains $\mathcal{H}$ as a subclass and
$\mathcal{H}$ is said to be \emph{closed} if $\overline{\mathcal{H}}=\mathcal{H}.$
A closed class $\mathcal{H}$ is said to be \emph{projective} if,
for each object $C\in\mathfrak{C}$, there is a morphism $f:P\rightarrow C$ in
$\mathcal{H}$ where $P$ is $\mathcal{H}$-projective.

If $\mathfrak{C}$ is an abelian category and
$\mathcal{H}$ is a \emph{closed} class of morphisms in $\mathfrak{C}$, then
a morphism $f\in\mathfrak{C}$ is called $\mathcal{H}$
\emph{-admissible} if in the canonical factorization
$
f=\mu\circ\xi,
$
where $\mu$ is a monomorphism and $\xi$ is an epimorphism, $\xi$ is an element of
$\mathcal{H}.$ An exact sequence in $\mathfrak{C}$ is called
$\mathcal{H}$\emph{-exact} if all its morphisms are $\mathcal{H}
$-admissible. Finally, an $\mathcal{H}$\emph{-projective resolution
}of an object $C\in\mathfrak{C}$ is an $\mathcal{H}$-exact sequence 
\[
\cdots\longrightarrow P_{n}\overset{d_{n}}{\longrightarrow}P_{n-1}
\overset{d_{n-1}}{\longrightarrow}\cdots\overset{d_{2}}{\longrightarrow}
P_{1}\overset{d_{1}}{\longrightarrow}P_{0}\overset{d_{0}}{\longrightarrow
}C\longrightarrow0 ,
\]
such that $P_{n}$ is $\mathcal{H}$-projective, for every
$n\in\mathbb{N}$. If $\mathcal{H}$ is a projective class of
epimorphisms in $\mathfrak{C}$, then
every object in $\mathfrak{C}$ admits an $\mathcal{H}$-projective resolution.

Let $\mathfrak{B},$ $\mathfrak{C}$ be abelian
categories and let $\mathcal{H}$ be a projective class of epimorphisms in
$\mathfrak{B}$ (so that every object in $\mathfrak{B}$ admits an
$\mathcal{H}$-projective resolution).  Given a contravariant
additive functor $\mathbf{T:}\mathfrak{B}\rightarrow\mathfrak{C}$ and given an
$\mathcal{H}$-projective resolution in $\mathfrak{B}$
\[
\mathbf{P}_{\bullet}\longrightarrow B\longrightarrow0
\]
of $B$, the object $\mathbf{H}^{n}(\mathbf{T}(\mathbf{P}_{\bullet}))$ depends only on
$B\ $and yields an additive functor
\[
\mathrm{R}_{\mathcal{H}}^{n}\mathbf{T}:\mathfrak{B}\rightarrow\mathfrak{C}
,\mathfrak{\qquad}\mathrm{R}_{\mathcal{H}}^{n}\mathbf{T}\left(  B\right)
:=\mathbf{H}^{n}(\mathbf{T}(\mathbf{P}_{\bullet})).
\]
The functor $\mathrm{R}_{\mathcal{H}}^{n}\mathbf{T}$ is called the $n$-th
\emph{right }$\mathcal{H}$\emph{-derived functor} of $\mathbf{T}.$

Similarly to a non-relative case, any short $\mathcal{H}$-exact sequence
\[
0\rightarrow B_{1}\rightarrow B_{2}\rightarrow B_{3}\rightarrow0
\]
in $\mathfrak{B}$ yields a long exact sequence
\begin{multline*}
0\rightarrow\mathrm{R}_{\mathcal{H}}^{0}\mathbf{T}\left(  B_{3}\right)
\rightarrow\mathrm{R}_{\mathcal{H}}^{0}\mathbf{T}\left(  B_{2}\right)
\rightarrow\mathrm{R}_{\mathcal{H}}^{0}\mathbf{T}\left(  B_{1}\right)
{\rightarrow}\mathrm{R}_{\mathcal{H}}^{1}\mathbf{T}\left(
B_{3}\right)  \rightarrow\cdots\\
\cdots{\rightarrow}\mathrm{R}_{\mathcal{H}}
^{n}\mathbf{T}\left(  B_{3}\right)  \rightarrow\mathrm{R}_{\mathcal{H}}
^{n}\mathbf{T}\left(  B_{2}\right)  \rightarrow\mathrm{R}_{\mathcal{H}}
^{n}\mathbf{T}\left(  B_{1}\right)  {\rightarrow
}\mathrm{R}_{\mathcal{H}}^{n+1}\mathbf{T}\left(  B_{3}\right)  \rightarrow
\cdots
\end{multline*}
of $\mathcal{H}$-derived functors (cf.\ \cite[Theorem 2.1, page 309]{HS}).

Let $\mathfrak{B},$ $\mathfrak{C}$ be abelian categories and let $\mathcal{H}
$ be a projective class of epimorphisms in $\mathfrak{B}$. A
contravariant functor $\mathbf{T}:\mathfrak{B}\rightarrow\mathfrak{C}$ is
said to be \emph{left}\textbf{\ }$\mathcal{H}$\emph{-exact}\textbf{\ }if, for
every $\mathcal{H}$\textbf{-}exact sequence
$
B_{1}\rightarrow B_{2}\rightarrow B_{3}\rightarrow0,
$
the sequence $0\rightarrow\mathbf{T}\left(  B_{3}\right)  \rightarrow
\mathbf{T}\left(  B_{2}\right)  \rightarrow\mathbf{T}\left(  B_{1}\right)  $
is exact.
By \cite[pages 311--312]{HS} a contravariant left $\mathcal{H}$-exact functor
 $\mathbf{T}:\mathfrak{B}\rightarrow
\mathfrak{C}$, is additive and naturally isomorphic to
$\mathrm{R}_{\mathcal{H}}^{0}\mathbf{T}$. Furthermore,
$\mathrm{R}_{\mathcal{H}}^{n}\mathbf{T}\left(  P\right)  =0,$ for every
$n>0 $ and for every $\mathcal{H}$-projective object $P.$

We now provide the main example of closed projective class we are
interested in.

\begin{theorem}
\cite[Theorem 2.2]{Ar1}\label{teo 4.2.33} Let $\mathbb{H}:\mathfrak{B}
\rightarrow\mathfrak{A}$ be a covariant functor and consider the  class of
$\mathbb{H}$-relatively split morphisms:
\[
\mathcal{E}_{\mathbb{H}}:=\{f\in\mathfrak{B}\mid\mathbb{H}(f)\text{ splits in
}\mathfrak{A}\}.
\]
Let $\mathbb{T}:\mathfrak{A}\rightarrow\mathfrak{B}$ be a left adjoint of
$\mathbb{H}$ and let $\varepsilon:\mathbb{TH}\rightarrow\mathrm{Id}
_{\mathfrak{B}}$ be the counit of the adjunction. Then, for any object
$P\in\mathfrak{B}$, the following assertions are equivalent:
\begin{blist}
\item $P$ is $\mathcal{E}_{\mathbb{H}}$-projective.
\item  Every morphism $f:B\rightarrow P$ in $\mathcal{E}_{\mathbb{H}}$ has a section.
\item The counit $\varepsilon_{P}:\mathbb{TH}(P)\rightarrow P$ has a section.
\item  There is a splitting morphism $\pi:\mathbb{T}(X)\rightarrow P$ for a
suitable object $X\in\mathfrak{A}$.
\end{blist}
 In particular all objects of the
form $\mathbb{T}(X)$, $X\in\mathfrak{A}$, are $\mathcal{E}_{\mathbb{H}}
$-projective. Moreover $\mathcal{E}_{\mathbb{H}}$ is a closed
projective class.
\end{theorem}

Thus $\mathcal{E}_{\mathbb{H}}$ is a projective
class. Note that since, for any object $Y\in\mathfrak{B}$, the morphism
$\mathbb{H}\left(  \varepsilon_{Y}\right)$ is split by
$\eta_{\mathbb{H}\left(
Y\right)  }$, the counit of adjunction
$\varepsilon_{Y}$ is in the class $\mathcal{E}_{\mathbb{H}}$. To apply the derived functors one needs to determine, when $\mathcal{E}_{\mathbb{H}}$ is a class of epimorphisms (in which case
any object in $\mathfrak{B}$ admits an
$\mathcal{E}_{\mathbb{H}}$-projective resolution). The necessary and sufficient
conditions for this are given in the next proposition, which is the only (mildly) new result
in this section.

\begin{proposition}
\label{pro: epiclass} Let $\left(  \mathbb{T},\mathbb{H}\right)  $ be an
adjunction, where $\mathbb{H}:\mathfrak{B}\rightarrow\mathfrak{A}$ is a
covariant functor. Let $\varepsilon:
\mathbb{TH}\rightarrow\mathrm{Id}_{\mathfrak{B}}$ be the
counit of the adjunction and let $\eta:\mathrm{Id}_{\mathfrak{A}}
\rightarrow\mathbb{HT}$ be the unit of the adjunction. Consider the class of $\mathbb{H}$-relatively split morphisms
$$
\mathcal{E}_{\mathbb{H}}    :=\{f\in\mathfrak{B}\mid\mathbb{H}(f)\text{
splits in }\mathfrak{A}\}.
$$
The following assertions are equivalent:

\begin{blist}
\item  $\mathcal{E}_{\mathbb{H}}$ is a class of epimorphisms.

\item The counit $\varepsilon_{Y}:\mathbb{T}\mathbb{H}\left(  Y\right)  \rightarrow
Y$ is an epimorphism for every $Y\in\mathfrak{B}$.

\item $\mathbb{H}:\mathfrak{B}\rightarrow\mathfrak{A}$ is faithful.
\end{blist}
\end{proposition}
\begin{proof}
(a) $\Rightarrow$ (b)
For all objects $Y\in\mathfrak{B}$,
$\varepsilon_{Y}\in\mathcal{E}_{\mathbb{H}}$. Since $\mathcal{E}_{\mathbb{H}}$
is assumed to be a class of epimorphisms,  $\varepsilon_{Y}$ is an
epimorphism.

(b) $\Leftrightarrow$ (c) This is a standard description of a right adjoint faithful
functor, see e.g.\
\cite[Section~2.12,~Proposition~3]{Par:cat}.

(c) $\Rightarrow$ (a) It follows by the fact that faithful functors reflect
epimorphisms.
\end{proof}

By Theorem \ref{teo 4.2.33}, $\mathcal{E}_{\mathbb{H}}$ is always a projective class, and it
is a class of epimorphisms, provided the equivalent conditions of
Proposition~\ref{pro: epiclass} hold. In this case any object in $\mathfrak{B}$
admits an
$\mathcal{E}_{\mathbb{H}}$-projective resolution which  is unique up to a homotopy.
Thus, for every $B^{\prime
}\in\mathfrak{B}$, one can consider the right
$\mathcal{E}_{\mathbb{H}}$-derived functors
$\mathrm{R}_{\mathcal{E}_{\mathbb{H}}}^{\bullet}\F_{B^{\prime}}$ of $\F_{B^{\prime
}}:=\mathfrak{B}(-,B^{\prime}):\mathfrak{B}\rightarrow\mathfrak{Ab}$. These functors
play a special role in what follows.

\begin{definition}
\label{def: Ext} Let $\mathfrak{A,B}$ be abelian categories. Let $\left(
\mathbb{T},\mathbb{H}\right)  $ be an adjunction, where $\mathbb{H}
:\mathfrak{B}\rightarrow\mathfrak{A}$ is a covariant functor. If
$\mathcal{E}_{\mathbb{H}}$ is a class of epimorphisms and the functor
$\F_{B^{\prime
}}:=\mathfrak{B}(-,B^{\prime})$ is left $\mathcal{E}_{\mathbb{H}}$-exact for every
$B^{\prime}\in\mathfrak{B}$, then for every $B,B^{\prime}\in\mathfrak{B}$, we set:
\[
\mathrm{Ext}_{\mathcal{E}_{\mathbb{H}}}^{\bullet}(B,B^{\prime})=\mathrm{R}
_{\mathcal{\mathcal{E}_{\mathbb{H}}}}^{\bullet}\F_{B^{\prime}}(B).
\]
\end{definition}

The  study of relative injectivity  can be carried in a dual way,
i.e.\  working in the opposite category of $\mathfrak{C}$ (note
that if a category is abelian, so is its opposite category).
In particular, the dual of
Theorem~\ref{teo 4.2.33}, \cite[Theorem 2.3]{Ar1},  states that  the
class of relatively cosplit morphisms,
\[
\mathcal{I}_{\mathbb{T}}:=\{g\in\mathfrak{A}\mid\mathbb{T}(g)\text{
cosplits in }\mathfrak{B}\},
\]
is a closed injective class. Dualizing Proposition~\ref{pro: epiclass} one
concludes that   $\mathcal{I}_{\mathbb{T}}$ is a class of monomorphisms iff $\mathbb{T}$
is a faithful functor.

\subsection{Separable functors}
The notion of a separable functor was introduced in \cite{NdO}.
Following the formulation in \cite{Rafael}, a covariant functor
 $\mathbb{H}:\mathfrak{B}\rightarrow\mathfrak{A}$ is said to be {\em separable}
 iff
the transformation
$
{\mathfrak{B}}(-,-)\rightarrow {\mathfrak{A}
}({\mathbb{H}}(-),{\mathbb{H}}(-)),$ $f\mapsto {\mathbb{H}}({f}),
$
 is a split natural monomorphism.

As explained in \cite[Lemma~1.1]{NdO}, any equivalence of categories is separable,
 and a composition of separable
functors is separable. Furthermore if a functor $\mathbb{H}\circ\mathbb{T}$ is
separable, then so is $\mathbb{T}$. By \cite[Proposition 1.2]{NdO}, a separable functor reflects split
monomorphisms and
split epimorphisms. This then implies that, for a pair of functors
$\mathbb{T}:\mathfrak{A}
\rightarrow\mathfrak{B}$ and $\mathbb{H}:\mathfrak{B}\rightarrow\mathfrak{C}$,
with $\mathbb{H}$ separable, the class of $\mathbb{H}\circ \mathbb{T}$-relatively
split morphisms  (resp.\ $\mathbb{H}\circ \mathbb{T}$-relatively
cosplit morphisms) is the same
as the class of $\mathbb{T}$-relatively
split morphisms (resp.\
$\mathbb{T}$-relatively cosplit morphisms), i.e.
\[
\mathcal{E}_{\mathbb{H}\circ \mathbb{T}} = \mathcal{E}_{\mathbb{T}} \qquad
\mbox{ (resp.} \;\mathcal{I}_{\mathbb{H}\circ \mathbb{T}} =
\mathcal{I}_{\mathbb{T}} \mbox{)}.
\]

A particularly useful criterion  of separability of a functor with
an adjoint is provided by the  Rafael Theorem:

\begin{theorem}
\label{teo Rafael}\cite[Theorem 1.2]{Rafael} Let $\mathbb{T}$ be a left adjoint
of a covariant functor $\mathbb{H}$.
\begin{zlist}
\item $\mathbb{T}$ is separable if and only if the unit of adjunction is a natural
section.
\item $\mathbb{H}$ is separable if and only if the counit of adjunction is a natural
retraction.
\end{zlist}
\end{theorem}

Combining Theorem~\ref{teo Rafael} with Theorem~\ref{teo 4.2.33} (and its dual)
we obtain
\begin{corollary}
\label{coro Rafael} Let
$\mathbb{T}:\mathfrak{A}\rightarrow\mathfrak{B}$ be a covariant functor with right
adjoint $\mathbb{H}$.
\begin{zlist}
\item If $\mathbb{H}$ separable, then any object in $\mathfrak{B}$ is
$\mathcal{E}_{\mathbb{H}}$-projective.
\item If
 $\mathbb{T}$ separable, then  any object in $\mathfrak{A}$ is
$\mathcal{I}_{\mathbb{T}}$-injective.
\end{zlist}
\end{corollary}

\section{Module-relative Hochschild cohomology}\label{sec.mod}
In this section we introduce and compute (in a special case)
the Hochschild cohomology relative to a bimodule. This cohomology is used in
the description of separable and formally smooth bimodules.

Let $A,$ $B$ and $T$ be rings. Given a bimodule ${_{B}M_{A}}$, consider the
following adjunction
\[
\begin{tabular}
[c]{lll}
$\mathbb{L}_{T}:{_{A}\mathcal{M}_{T}}\rightarrow{_{B}\mathcal{M}_{T},}$ &  &
${\mathbb{R}}_{T}:{_{B}\mathcal{M}_{T}}\rightarrow{_{A}\mathcal{M}_{T}}$\\
$\mathbb{L}_{T}\left(  X\right)  =M\otimes_{A}X{,}$ &  & $\mathbb{R}
_{T}\left(  Y\right)  ={_{B}\mathrm{Hom}}\left(  M,Y\right)
{,}$
\end{tabular}
\ \ \ \
\]
Note that the counit $\varepsilon^T$ of this adjunction is, for all
$Y\in {}_{B}\mathcal{M}_{T}$,
$$
\varepsilon^T_Y : M\otimes_A {}_{B}\mathrm{Hom}\left(  M,Y\right)
\to Y, \qquad m\otimes_A f\mapsto (m)\! f.
$$
We would like to compute the cohomology relative to the class
$$
\eclass {M,T} := \mathcal{E}_{\mathbb{R}_T} = \{f \in  {}_{B}\mathcal{M}_{T}\;
|\; {}_{B}\mathrm{Hom}\left( M, f\right) \mbox{ splits in
}{}_{A}\mathcal{M}_{T}\}
$$
To apply the derived functors we need to determine, when $\eclass T$ is a class
of epimorphisms.

\begin{proposition}
\label{pro: generator}
Let $\varepsilon^{T}:\mathbb{L}_{T}{\mathbb{R}}
_{T}\rightarrow\mathrm{Id}_{{_{B}\mathcal{M}_{T}}}$ be the counit of the
adjunction $\left(  \mathbb{L}_{T},{\mathbb{R}}_{T}\right)  $.
The following assertions are equivalent:

\begin{enumerate}
\item[\em (a)] $\mathcal{E}_{{M,T}}$ is a class of epimorphisms for
every ring $T$.

\item[\em (a$'$)] $\mathcal{E}_{{M,B}}$ is a class of epimorphisms.

\item[\em (a$''$)] $\mathcal{E}_{M,{\mathbb{Z}}}$ is a
class of epimorphisms.

\item[\em (b)] The counit $\varepsilon_{Y}^{T}:\mathbb{L}_{T}{\mathbb{R}}_{T}\left(
Y\right)  \rightarrow Y$ is an epimorphism for every ring $T$ and for every
$Y\in{_{B}\mathcal{M}_{T}}$.

\item[\em (b$'$)] The counit $\varepsilon_{Y}^{B}:\mathbb{L}_{B}{\mathbb{R}}_{B}\left(
Y\right)  \rightarrow Y$ is an epimorphism for every
$Y\in{_{B}\mathcal{M}_{B}}$.

\item[{\em (b$''$)}] The counit $\varepsilon_{Y}^{\mathbb{Z}}:\mathbb{L}
_{\mathbb{Z}}{\mathbb{R}}_{\mathbb{Z}}\left(  Y\right)  \rightarrow Y$ is an
epimorphism for every $Y\in{_{B}\mathcal{M}_{\mathbb{Z}}} =
{_{B}\mathcal{M}}$.

\item[\em (c)] ${\mathbb{R}}_{T}:{_{B}\mathcal{M}_{T}}\rightarrow{_{A}
\mathcal{M}_{T}}$ is faithful for every ring $T$.

\item[\em (c$'$)] ${\mathbb{R}}_{B}:{_{B}\mathcal{M}_{B}}\rightarrow{_{A}
\mathcal{M}_{B}}$ is faithful.

\item[\em (c$''$)] ${\mathbb{R}}_{\mathbb{Z}}:{_{B}\mathcal{M}
}\rightarrow{_{A}\mathcal{M}}$ is faithful.

\item[\em (d)] The evaluation map
\[
\mathrm{ev}_{M}:M\otimes_{A}{^{\ast}\!M}\rightarrow B,\qquad\mathrm{ev}
_{M}\left(  m\otimes_{A}f\right)  =\left(  m\right) \! f,
\]
where
${}{^{\ast}\!M}:={{}_{B}\mathrm{Hom}}\left( M,B\right)$
is an epimorphism (of $B$-bimodules).

\item[\em (e)] $M$ is a generator in ${_{B}\mathcal{M}}$.
\end{enumerate}
\end{proposition}

\begin{proof}
The equivalences
$\left(
a\right)  \Leftrightarrow \left(  b\right)  \Leftrightarrow \left(  c\right)  $,
$\left(  a^{\prime
}\right) \Leftrightarrow  \left(  b^{\prime}\right)  \Leftrightarrow \left(  c^{\prime}\right)$
and
$\left(  a^{\prime\prime}\right) \Leftrightarrow \left(  b^{\prime\prime}\right)
 \Leftrightarrow\left(  c^{\prime\prime}\right)  $ follow by
 Proposition~\ref{pro: epiclass}. The implication
 $\left(  c\right)  \Rightarrow\left(  c^{\prime}\right)  $  is
obvious, while $\left(  a^{\prime}\right)  \Rightarrow\left(  d\right)$ follows by
identifying $\mathrm{ev}_{M}$ with the counit of adjunction (at $B$) $\varepsilon_{B}^{B}\in\mathcal{E}_{M,{B}}$.
The latter is in the class $\eclass {M,B}$, hence is an epimorphism
(by assumption $\left(  a^{\prime}\right)$).  The equivalences
$\left(
c^{\prime\prime}\right)  \Leftrightarrow
\left(  d\right)  \Leftrightarrow\left(  e\right)$ are standard
characterizations of generators in the category of modules
(cf.\ \cite[13.7]{Wisbauer-book}).
Finally, since, for
all $f\in{_{B}\mathrm{Hom}}_{T}\left(  Y,Y^{\prime}\right)
$,  $\mathbb{R}_{\mathbb{Z}}\left(  f\right)  ={_{B}\mathrm{Hom}}\left(
M,f\right)  ={\mathbb{R}}_{T}\left(  f\right)$, the condition
$\left(  c^{\prime\prime}\right)$ implies $\left(  c\right)  $.
\end{proof}

Clearly, for every $Y^{\prime}\in{_{B}\mathcal{M}_{B},}$ the functor
$\F_{Y^{\prime}}:={_{B}\mathrm{Hom}}_{B}(-,Y^{\prime}):{_{B}\mathcal{M}_{B}
}\rightarrow\mathfrak{Ab}$ is left $\mathcal{E}_{M,{B}}$-exact so, in view
of equivalent conditions in Proposition~\ref{pro: generator} we can propose the following
\begin{definition}
\label{def.Hoch}
Let $M$ be a $B$-$A$ bimodule which is a generator as a left $B$-module, and
let $\eclass {M, B}$ be the class of all $B$-bimodule maps $f$, such that ${}_{B}\mathrm{Hom}\left( M, f\right)$ splits as an $A$-$B$ bimodule map.
The \emph{$M$-Hochschild
cohomology of $B$ with coefficients in a $B$-bimodule $N$} is
defined by
\[
\mathrm{H}_{M}^{\bullet}(B,N):=\mathrm{Ext}_{\mathcal{E}_{M,{B}}}
^{\bullet}(B,N),
\]
(cf.\ Definition~\ref{def: Ext} for the explanation of the relative Ext-functor).

If the number
\[
 \min\left\{
n\in\mathbb{N}\cup \{0\}\mid\mathrm{H}_{M}^{n+1}\left(  B,N\right)
=0\text{ for every }N\in{_{B}}\mathcal{M}{_{B}}\right\}
\]
exists, then it is called an \emph{$M$
-Hochschild dimension of $B$} and is denoted by $\mathrm{Hdim}_{M}\left(  B\right)$.
Otherwise $B$ is said to have an infinite $M$-Hochschild dimension.
\end{definition}

Similarly to the non-relative case, $M$-Hochschild cohomology can be equivalently
described as the cohomology of a complex associated to the standard resolution.
The standard resolution can be described in general as follows.
Start with an additive functor
$\mathbb{H}:\mathfrak{B}\rightarrow\mathfrak{A}$ of abelian categories
with a left adjoint $\mathbb{T}$.
This defines a comonad $F := \mathbb{TH}$ on
 $\mathfrak{B}$ with the counit given by the counit of adjunction
 $(\mathbb{T},\mathbb{H})$,
$\varepsilon:\mathbb{TH}\rightarrow\mathrm{Id}_{\mathfrak{B}}$.
For an object $B\in\mathfrak{B}$, one considers the
 \emph{associated augmented chain
complex}
\[
\cdots\overset{d_{3}}{\longrightarrow}F^{3}(B)\overset{d_{2}}{\longrightarrow
}F^{2}(B)\overset{d_{1}}{\longrightarrow}F(B)\overset{d_{0}}{\longrightarrow}
F^{0}(B):=B\rightarrow0
\]
where
\[
d_{n}=
{\displaystyle\sum\limits_{i=0}^{n}}
\left(  -1\right)  ^{i}F^{i}\left(  \varepsilon_{F^{n-i}(B)}\right)
\]
(see \cite[8.6.4, page 280]{Weibel}).

\begin{proposition}
Let $\mathfrak{A,B}$ be abelian categories. Let
$\mathbb{H}:\mathfrak{B}
\rightarrow\mathfrak{A}$ be a faithful covariant functor
with a left adjoint $\mathbb{T}$.
For all objects of $\mathfrak{B}$,
the associated augmented chain complex
is an $\mathcal{E}_{\mathbb{H}}$-exact sequence.
\end{proposition}

\begin{proof}
Let $\varepsilon:\mathbb{TH}\rightarrow\mathrm{Id}_{\mathfrak{B}}$ be the
counit of the adjunction and let $\eta:\mathrm{Id}_{\mathfrak{A}}
\rightarrow\mathbb{HT}$ be the unit of the adjunction. For all integers $n\geq -1$,
define
\[
s_{n}:=\eta_{\mathbb{H}F^{n+1}(B)}:\mathbb{H}F^{n+1}(B)\rightarrow\mathbb{H}
F^{n+2}(B).
\]
Then,
$
\mathbb{H}\left(  d_{0}\right)  \circ s_{-1}=\mathbb{H}\left(  \varepsilon
_{B}\right)  \circ\eta_{\mathbb{H}(B)}=\mathrm{Id}_{\mathbb{H}(B)}\text{.}
$
Furthermore, for all $n\geq0$,
$d_{n}  = \varepsilon_{F^{n}(B)}-F\left(  d_{n-1}\right)$,
so that
\begin{align*}
\mathbb{H}\left(  d_{n+1}\right)  \circ s_{n}
  & =\mathbb{H}\left(   \varepsilon_{F^{n+1}(B)}\right)  \circ\eta_{\mathbb{H}
F^{n+1}(B)}-\mathbb{H}F\left(  d_{n}\right)  \circ\eta_{\mathbb{H}F^{n+1}(B)}\\
&  =\mathbb{H}\left(  \varepsilon_{F^{n+1}(B)}\right)  \circ\eta_{\mathbb{H}
F^{n+1}(B)}-\eta_{\mathbb{H}F^{n}(B)}\circ\mathbb{H}\left(  d_{n}\right) \\
&  =\mathrm{Id}_{\mathbb{H}F^{n+1}(B)}-s_{n-1}\circ\mathbb{H}\left(
d_{n}\right)  ,
\end{align*}
where the second equality follows by the naturality of the unit of adjunction.
Hence $s_\bullet$ is a contracting homotopy for the complex
$(\mathbb{H}(F^\bullet (B)), \mathbb{H}(d_\bullet))$,
which implies that the augmented chain complex
$(F^\bullet (B), d_\bullet)$
is an $\mathcal{E}_{\mathbb{H}}$-exact sequence.
\end{proof}

In the case of the adjunction $(\mathbb{L}_B,\mathbb{R}_B)$, the comonad
is $F=M\otimes_A {}_B\mathrm{Hom} (M,-)$. Application of  the functor ${_{B}\mathrm{Hom}}_{B}(-,N):{_{B}\mathcal{M}_{B}
}\rightarrow\mathfrak{Ab}$ to the associated augmented chain complex,
results in the cochain complex
$$
({}_{B}{\mathrm{Hom}}_{B}(F^\bullet (B),N),
d^{\bullet}:={}_{B}\mathrm{Hom}_{B}(d_{\bullet},N)),
$$
whose cohomology is $\mathrm{H}^\bullet_M(B,N)$.

The $M$-Hochschild cohomology has a particularly simple description in the case
$M$ is a progenerator left $B$-module. In this case it can be identified with a (relative)
Hochschild cohomology of the endomorphism ring of $M$. This can be
 described as follows.

Given a ring extension $A\to S$ (or an $A$-ring $S$), the $A$-relative
Hochschild cohomology
of $S$ with values in an $S$-bimodule $W$ \cite{Hoc:rel},
$\mathrm{H}^\bullet(S|A, W)$, is defined as the
cohomology of the  cochain complex
\[
0\rightarrow \!{_{A}\mathrm{Hom}}_{A}(A,\! W)\!\overset{b^{0}}{\rightarrow}\!
{_{A}\mathrm{Hom}}_{A}({S},\! W)\!\overset{b^{1}}{\rightarrow}\!{_{A}\mathrm{Hom}
}_{A}(S^{\otimes_{A}2},\! W)\!\overset{b^{2}}{\rightarrow}\!{_{A}\mathrm{Hom}}
_{A}(S^{\otimes_{A}3},\! W)\!\overset{b^{3}}{\rightarrow}\!\cdots ,
\]
where, for all $f\in  {}_{A}\mathrm{Hom}_{A}({S}^{\otimes_A n},W)$,
$n=0,1,2,\ldots$,
$$
b^n(f) = \mu_{W}^{l}\circ\left(  S\otimes_{A}f\right)  +
\sum_{i=1}^n(-1)^i f\circ(S^{\otimes i-1}\otimes_{A}m_{S}\otimes_{A}S^{\otimes n-i})
 + (-1)^{n+1}
\mu_{W}^{r}\circ(f\otimes_{A}S).
$$
Here $\mu_W^l$, $\mu_W^r$ denote left, respectively, right $S$-multiplication on
$W$ and $m_S: S\otimes_A S\to S$ is the product map. Also, in case $n=0$,
the obvious isomorphisms $A\otimes_AS\simeq S\otimes_AA\simeq S$ are implicitly
used. $\mathrm{H}^\bullet(S|A, W)$ can be understood as the Hochschild
cohomology of the
algebra $S$ in monoidal category of $A$-bimodules
(cf.\ \cite[Theorem~4.42]{AMS1}). The {\em
Hochschild dimension of
 $S$ over $A$}  is then defined by
\[
\mathrm{Hdim}\left(  S | A\right) := \min\left\{  n\in\mathbb{N}\cup \{0\}\mid\mathrm{H}^{n+1}\left(  S|A,W\right)  =0\text{
for every }W\in{_{S}}\mathcal{M}{_{S}}\right\} ,
\]
provided that the minimum on the right hand side exists.
\begin{theorem}
\label{teo: Hoch S} Let $A,$ $B$ be rings. Consider a bimodule ${_{B}M_{A}}$
such that ${_{B}M}$ is a progenerator. Let $S$ be the endomorphism ring
of the left $B$-module $M$.
Then, for all  $B$-bimodules $N$,
\[
\mathrm{H}_{M}^{\bullet}(B,N)=\mathrm{H}^{\bullet}\left(  S|A,{^{\ast}}M\otimes
_{B}N\otimes_{B}M\right) .
\]

Furthermore, for a fixed $n\in\mathbb{N}$, the following assertions are
equivalent:
\begin{blist}
\item $\mathrm{H}_{M}^{n}(B,N)=0,$ for every $B$-bimodule $N$.

\item $\mathrm{H}^{n}\left(  S|A,{W}\right)  =0,$ for every $S$-bimodule $W$.
\end{blist}

In particular
\[
\mathrm{Hdim}_{M}\left(  B\right)  =\mathrm{Hdim}\left(
S|A\right)  .
\]
\end{theorem}
\begin{proof}
Since $M$ is a finitely generated and projective left $B$-module, the functor
${\mathbb{R}}_{B}$ is isomorphic to  ${^{\ast}M}\otimes_{B}\left(
-\right)  _{B}.$
The comonad comes out as
\[
F\left(  N\right)  =\mathbb{L}_{B}\mathbb{R}_{B}\left(  N\right)
\equiv {C}\otimes_{B}N\text{, for all }N\in{_{B}\mathcal{M}_{B}},
\]
where
${C}:=M\otimes_{A}{^{\ast}}M\equiv\mathbb{L}_{B}\mathbb{R}_{B}\left(
B\right) $. The counit of adjunction at $B$ is simply the evaluation map
$\mathrm{ev}_{M}:M\otimes_{A}{^{\ast}\!M}\rightarrow B$,
$ m\otimes_{A}f \mapsto \left(  m\right) \! f$.
Using the standard isomorphisms $C\otimes_B B \simeq C$, and
applying the Hom-functor to
the augmented chain complex associated to $B$, we can identify
$\mathrm{H}^\bullet_M(B,N)$ with the cohomology of the cochain complex
\[
0\rightarrow \!{_{B}\mathrm{Hom}}_{B}(B,\! N)\!\overset{d^*_{0}}{\rightarrow}\!
{_{B}\mathrm{Hom}}_{B}({C},\! N)\!\overset{d^*_{1}}{\rightarrow}\!{_{B}\mathrm{Hom}
}_{B}(C^{\otimes_{B}2},\! N)\!\overset{d^*_{2}}{\rightarrow}\!{_{B}\mathrm{Hom}}
_{B}(C^{\otimes_{B}3},\! N)\!\overset{d^*_{3}}{\rightarrow}\!\cdots ,
\]
where, for all $f\in  {}_{B}\mathrm{Hom}_{B}({C}^{\otimes_B n},N)$,
$n=0,1,2,\ldots$,
\[
d_{n}^*(f)=
{\displaystyle\sum\limits_{i=0}^{n}}
\left(  -1\right)  ^{i} f\circ(C^{\otimes_ B i}\otimes_B \mathrm{ev}_{M}\otimes_B
  C^{\otimes_ Bn- i}).
\]
Since $M$ is a finitely generated and projective left $B$-module, $S$ can be
identified with ${}^*M\otimes_B M$. Under this identification, the product is given
by ${}^*M\otimes_B\mathrm{ev}_M\otimes_B M$ and the unit is the dual basis
element $\sum_{a\in I}{^{\ast}e_{a}}\otimes_{B}e_{a}$. Furthermore, one can
consider the isomorphisms
\[
\Phi_{n}:{_{B}\mathrm{Hom}}_{B}({C}^{\otimes_{B}n+1},N)\rightarrow
{_{A}\mathrm{Hom}}_{A}(S^{\otimes_{A}n},{^{\ast}}M\otimes_{B}N\otimes_{B}M),
\]
defined by
\[
\left[  \Phi_{n}\left(  f\right)  \right]  \left(  x\right)  =\left(  {^{\ast
}}M\otimes_{B}f\otimes_{B}M\right)  \left(  1_{S}\otimes_{A}x\otimes_{A}
1_{S}\right)  ,\text{ for every }x\in S^{\otimes_{A}n}.
\]
Using the definitions of cochain maps and above identification of $S$, one
easily checks that these isomorphisms fit into the commutative diagrams
\[
\begin{tabular}
[c]{ccc}
${_{B}\mathrm{Hom}}_{B}(B,N)$ & $\overset{d_{0}^{\ast}}{\longrightarrow}$ &
${_{B}\mathrm{Hom}}_{B}({C},N)$\\
$\simeq\downarrow$ &  & $\downarrow\Phi_{0}$\\
${N}^{B}=\left\{  n\in N\mid bn=nb,\forall b\in B\right\}  $ & $\overset
{b^{-1}}{\longrightarrow}$ & ${_{A}\mathrm{Hom}}_{A}(A,{^{\ast}}M\otimes
_{B}N\otimes_{B}M)$
\end{tabular}
\ \ \ \
\]
and
\[
\begin{tabular}
[c]{ccc}
${_{B}\mathrm{Hom}}_{B}({C}^{\otimes_{B}n},N)$ & $\overset
{d_{n}^{\ast}}{\longrightarrow}$ & ${_{B}\mathrm{Hom}}_{B}({C}
^{\otimes_{B}n+1},N)$\\
$\Phi_{n-1}\downarrow$ &  & $\downarrow\Phi_{n}$\\
${_{A}\mathrm{Hom}}_{A}(S^{\otimes_{A}n-1},{^{\ast}}M\otimes_{B}N\otimes
_{B}M)$ & $\overset{b^{n-1}}{\longrightarrow}$ & ${_{A}\mathrm{Hom}}
_{A}(S^{\otimes_{A}n},{^{\ast}}M\otimes_{B}N\otimes_{B}M)$
\end{tabular}
\ \ \ \
\]
This immediately implies that
\[
\mathrm{H}_{M}^{\bullet}(B,N)=\mathrm{H}^{\bullet}\left(  S|A,{^{\ast}}M\otimes
_{B}N\otimes_{B}M\right) ,
\]
as required.

 It remains to prove that the statements $\left(  a\right)  $ and
$\left(  b\right)  $ are equivalent. The implication $\left(  b\right)  \Rightarrow
\left(  a\right)  $ is obvious. To prove the converse,
take any $S$-bimodule $W$ and define a $B$-bimodule $N=M\otimes_{S}
W\otimes_{S}{^{\ast}}M.$ Then
\[
{^{\ast}}M\otimes_{B}N\otimes_{B}M={^{\ast}}M\otimes_{B}M\otimes_{S}
W\otimes_{S}{^{\ast}}M\otimes_{B}M=S\otimes_{S}W\otimes_{S}S\simeq W.
\]
This completes the proof.
\end{proof}

\section{Separable bimodules}
The aim of the section is to supplement (and extend) the functorial description
of separable bimodules in \cite[Corollary 5.8]{CK} by the cohomological
description of such bimodules. First recall from \cite{Su} the following

\begin{definition}
Let $A,$ $B$ be rings. A $B$-$A$ bimodule $M$
 is said to be \emph{separable}, or $B$ is said to be
{\em $M$-separable over $A$} if the evaluation map
\[
\mathrm{ev}_{M}:M\otimes_{A}{^{\ast}M}\rightarrow B,\qquad\mathrm{ev}
_{M}\left(  m\otimes_{A}f\right)  =\left(  m\right)\!  f,
\]
is a split epimorphism
of $B$-bimodules.
\end{definition}
Throughout this section, $M$ is a $B$-$A$ bimodule,  and
$\mathbb{L}_T$, $\mathbb{R}_T$, $\eclass {M,T}$
are the functors and the class of morphisms (associated to $M$)
described at the beginning
of Section~\ref{sec.mod}.

\begin{proposition}
\label{pro: SepBimod} The following assertions are
equivalent for a $B$-$A$ bimodule $M.$

\begin{blist}
\item $M$ is a separable bimodule.

\item For all rings $T$, ${\mathbb{R}}_{T}:{_{B}\mathcal{M}_{T}}\rightarrow{_{A}\mathcal{M}_{T}}$
is a separable functor.

\item ${\mathbb{R}}_{B}:{_{B}\mathcal{M}_{B}}\rightarrow{_{A}\mathcal{M}_{B}}$
is a separable functor.

\item ${\mathbb{R}}_{\mathbb{Z}}:{_{B}\mathcal{M}}\rightarrow{_{A}\mathcal{M}
}$ is a separable functor.

\item Any $B$-bimodule is $\eclass {M,B}$-projective.

\item $B$ is $\eclass{M,B}$-projective.

\item $M$ is a generator in ${_{B}\mathcal{M}}$ and $\mathrm{H}_{M}
^{n}(B,N)=0,$ for every $N\in{_{B}\mathcal{M}_{B}}$ and for every $n\geq1.$

\item $M$ is a generator in ${_{B}\mathcal{M}}$ and $\mathrm{H}_{M}
^{1}(B,N)=0,$ for every $N\in{_{B}\mathcal{M}_{B}.}$

\item $M$ is a generator in ${_{B}\mathcal{M}}$ and
$\mathrm{Hdim}_M(B) =0$.
\end{blist}
\end{proposition}

\begin{proof}
(a) $\Rightarrow$  (b) For any $B$-$T$ bimodule $Y$, there is a (natural in $Y$)
$A$-$T$ bimodule map
$$
\hat{\xi}: {}^*M\otimes _B Y \to {}_B\mathrm{Hom}(M,Y), \quad
{}^*m\otimes_B y \mapsto [m\mapsto (m){}^*m y],
$$
(cf.\ \cite[Proposition~20.10]{AndFul:ring}). Tensoring this map with
$M$,  we obtain a $B$-$T$ bimodule
map
$$
\xi :  M\otimes_{A}{^{\ast}}M \ot_B Y\to M\otimes_A
 {}_B\mathrm{Hom}(M,Y).
 $$
It is clear from the definition and  naturality  of $\hat{\xi}$ that,
 for all $x\in \left(  M\otimes_{A}\!{^{\ast}}\!M\right)^{B} \!:=
\! \{ x\in  M\otimes_{A}\!{^{\ast}}\!M \mid
\forall b\in B,\; xb=bx\}$, the map
 $\xi(x\otimes_B -): Y \to  M\otimes_A {}_B\mathrm{Hom}(M,Y)$ is natural in $Y$ and
 $B$-$T$ bilinear.

If $M$ is a separable bimodule, then
 there exists $s\in\left(  M\otimes_{A}{^{\ast}}M\right)  ^{B}$ such that
 $\mathrm{ev}_M(s) =1_{B}$. One easily checks that $\xi(s\otimes_A-)$
is a natural splitting of the counit of the adjunction $\left(
\mathbb{L}_{T},{\mathbb{R}}_{T}\right)$.
Hence, by Rafael's Theorem~\ref{teo Rafael},
${\mathbb{R}}_{T}$ is a separable functor.

Implications
(b) $ \Rightarrow$ (c) and (e) $\Rightarrow$ (f) are obvious, while the equivalence
(a) $\Leftrightarrow$ (d)  is proven in \cite[Corollary 5.8]{CK}. The implication (c)$
 \Rightarrow$ (e) follows by Corollary~ \ref{coro Rafael}.

 Since $\mathrm{ev}_M$ is the same as the counit of adjunction
 $\left(
\mathbb{L}_{B},{\mathbb{R}}_{B}\right)$ evaluated at $B$, the
implication (f) $\Rightarrow$ (a) follows by Theorem~\ref{teo 4.2.33}.
Thus, if $B$ is $\eclass {M,B}$-projective, then
$\mathrm{ev}_{M}$ is an epimorphism, hence $M$ is a generator in ${_{B}\mathcal{M}}$ by Proposition~\ref{pro: generator}. Therefore, the
equivalences between  (f), (g), (h) and (i) follow by the definitions of $M$-Hochschild cohomology of $B$ with
coefficients in $N$ and $M$-Hochschild dimension of $B$, and by the properties of $\mathrm{Ext}_{\mathcal{E}
_{{M,B}}}^{\bullet}\left(  -,-\right)  .$
\end{proof}

Recall that a ring morphism $A\to S$ is called a {\em separable extension} if the
product map $m_S: S\otimes_AS\to S$ has an $S$-bimodule section.
\begin{proposition}
\cite[Theorem 1]{Su} \label{pro: sep S}Let $A,$ $B$ be rings. Consider a
bimodule ${_{B}M_{A}}$ such that ${M}$ is a finitely generated and projective
left $B$-module. Let
\[
S={_{B}\mathrm{End}}\left( M\right)  =\mathbb{R}
_{A}\left(  M\right)  \simeq{^{\ast}M\otimes}_{B}M.
\]
The following assertions are equivalent:
\begin{blist}
\item ${_{B}M_{A}}$ is a separable bimodule.

\item $M$ is a generator in ${_{B}\mathcal{M}}$ and the canonical morphism
$i:A\rightarrow S$, $a\mapsto [m\mapsto ma]$ is s a separable extension.
\end{blist}
\end{proposition}

\begin{proof}
 By Proposition~\ref{pro: SepBimod}, a separable bimodule $M$ is a generator in
${_{B}\mathcal{M}}$. Thus, in both cases, $M$ is a progenerator and, by
Theorem~\ref{teo: Hoch S},
$\mathrm{Hdim}_{M}\left(  B\right)  =\mathrm{Hdim}\left(
S|A\right)$. Since the extension $A\to S$ is separable if and only if
$\mathrm{Hdim}\left(
S|A\right) =0$ (cf.\ \cite[Theorem~4.43]{AMS1}), the assertion follows by Proposition~\ref{pro: SepBimod}.

The proposition can also be proven directly as follows. In both cases $M$ is
a progenerator left $B$-module, hence $S$ is isomorphic to ${^{\ast}}M\otimes_{B}M$ (${}_BM$ is finitely generated and projective) and $B$ is isomorphic to
$M\otimes_S{}^*\!M$ (${}_BM$ is a generator). These isomorphisms
allow one to identify the product $m_S$ in $S$ with
${^{\ast}}M\otimes_B\mathrm{ev}_M\otimes_{B}M$, and
$\mathrm{ev}_M$ with $M\otimes_Sm_S\otimes_S{}^*\!M$. With these
identifications, the mutual equivalence of statements (a) and (b) is clear.
\end{proof}

\begin{remark}
For a left $B$-module $M$, let
$S={{}_{B}\mathrm{End}}\left( M\right)$. For the $B$-$S$ bimodule
$M$, the map $i$ of Proposition~\ref{pro: sep S} is the identity and
hence it trivially defines a separable extension. Still ${_{B}M_{S}}$ needs
not to be a separable bimodule unless $_{B}M$ is a generator in ${_{B}
\mathcal{M}}$ (see Corollary~\ref{coro: M-static}).
\end{remark}

\section{Formally smooth bimodules}
In this section we introduce the notion of a formally smooth bimodule, and give
cohomological interpretation and describe examples of such bimodules.
Throughout this section, $M$ is a $B$-$A$ bimodule,  and, for any ring $T$,
$\mathbb{L}_T$, $\mathbb{R}_T$, $\eclass {M, T}$
are the functors and the class of morphisms (associated to $M$) described at the beginning
of Section~\ref{sec.mod}.

\begin{definition}
Let $A,$ $B$ be rings. A $B$-$A$ bimodule $M$ is said to be \emph{formally
smooth} or $B$ is said to be \emph{$M$-smooth over $A$}
whenever the kernel of the evaluation map 
\[
\mathrm{ev}_{M}:M\otimes_{A}{^{\ast}M}\rightarrow B,\qquad\mathrm{ev}
_{M}\left(  m\otimes_{A}f\right)  =\left(  m\right) \! f.
\]
is an
$\mathcal{E}_{M,B}$-projective $B$-bimodule.
\end{definition}

Following \cite{KonRos:non} a pair of functors
$\mathbb{U}_{\ast}:\bar{\mathfrak{A}}\rightarrow\mathfrak{A}$, $\mathbb{U}
^{\ast}:\mathfrak{A}\rightarrow\bar{\mathfrak{A}}$ such that $\mathbb{U}
^{\ast}$ is fully faithful and left adjoint to $\mathbb{U}_{\ast}$ is called a
\emph{Q-category}. As explained in \cite[Section~2.5]{KonRos:non}, to any
category $\mathfrak{C}$ and any class of morphisms $\mathcal{H}$ in
$\mathfrak{C}$ which contains all the identity morphisms, one can associate a
Q-category as follows. First construct the category $\mathfrak{H}$, whose
objects are elements $f$, $g$ of $\mathcal{H}$ and morphisms are commutative
squares
\[
\diagQcatA
\]
where the vertical arrows are  in $\mathfrak{C}$. The {\em direct image
functor} $\mathbb{U}^{\ast}$ is
\[
\mathbb{U}^{\ast}:M\mapsto\mathrm{id}_{M},\qquad \left(
\diagQcatB \right)  \mapsto\left(\vcenter{\diagQcatC}
\right)  .
\]
The {\em inverse image functor} $\mathbb{U}_{\ast}$ is defined by
\[
\mathbb{U}_{\ast}:\left(
\diagQcatB
\right)  \mapsto M,\qquad
\vcenter{\diagQcatA} \mapsto\left(\diagQcatB \right)  .
\]
We denote this Q-category by $\mathfrak{A}_{\mathcal{H}}$ and call it a
\emph{Q-category induced by the class $\mathcal{H}$}. Following
\cite[Sections~3.7 \& 4.5]{KonRos:non} an object $P\in\mathfrak{C}$ is said
to be
\emph{formally $\mathfrak{A}_{\mathcal{H}}$-smooth} if and only if, for every
$f\in\mathcal{H}$, the mapping $\mathfrak{C}(P,f)$ is a strict epimorphism
(i.e.\ a surjective map) of sets. Thus $P$ is formally $\mathfrak{A}
_{\mathcal{H}}$-smooth if and only if $P$ is $\mathcal{H}$-projective.
This leads immediately to the following lemma, which also explains the choice of the
terminology.

\begin{lemma}
A bimodule ${_{B}M_{A}}$ is {formally smooth} if and only if
$\mathrm{Ker}\left(  \mathrm{ev}_{M}\right)  $ is a formally smooth object
in the $Q$-category induced by the class of morphisms $\mathcal{E}_{M,B}$.
\end{lemma}

The following proposition gives the first examples of formally smooth bimodules.

\begin{proposition}\label{prop.1ex}
A $B$-$A$ bimodule $M$
is  formally smooth whenever one of the following conditions holds:
\begin{enumerate}
\item[(1)] $M$ is a separable bimodule.
\item[(2)] The map $\mathrm{ev}_{M}$ is injective.
\end{enumerate}
\end{proposition}

\begin{proof}
(1) If $M$ is a separable bimodule, then
Proposition~\ref{pro: SepBimod} implies that any $B$-bimodule is
$\eclass {M,B}$-projective. In particular $\mathrm{Ker}\left(
\mathrm{ev}_{M}\right)  $ is $\eclass {M,B}$-projective.

(2) If $\mathrm{ev}_{M}$ is injective, then $\mathrm{Ker}\left(
\mathrm{ev}_{M}\right)  $ is trivial, hence $\eclass {M,B}$-projective.
\end{proof}

The cohomological interpretation of formally smooth bimodules is provided by
the following
\begin{proposition}\label{prop.smooth}
Let $A,$ $B$ be rings. Take a $B$-$A$ bimodule $M$ which is a generator in ${_{B}\mathcal{M}}
$. Then the following assertions are equivalent:

\begin{blist}
\item $M$ is a formally smooth bimodule.
\item $\mathrm{H}_{M}^{n}(B,N)=0,$ for all $N\in{_{B}\mathcal{M}_{B}}$ and
all $n\geq2.$
\item $\mathrm{H}_{M}^{2}(B,N)=0,$ for all $N\in{_{B}\mathcal{M}_{B}.}$
\item $\mathrm{Hdim}_M(B) \leq 1$.
\end{blist}
\end{proposition}

\begin{proof}
The equivalences (b) $ \Leftrightarrow$ (c) $ \Leftrightarrow$ (d)
follow
 by the definitions
of the $M$-Hochschild cohomology of $B$ with coefficients in $N$ and
the $M$-Hochschild dimension of $B$,  and by the properties
of $\mathrm{Ext}_{\eclass{M,B}}^{\bullet}\left(  -,-\right).$

(a) $ \Leftrightarrow$ (c). Write $\left(
L,i_L\right)$ for the kernel of   $ \mathrm{ev}_{M}$,  and
consider  the following exact sequence of $B$-bimodules:
\[
\xymatrix{
0\ar[r] & L \ar[r]^{i_L\hspace{20pt}} & M\otimes_{A}{^{\ast}}\! M \ar[r]^{
\hspace{15pt}\mathrm{ev}_{M}} & B\ar[r]& 0.}
\]
Note that $\mathrm{ev}_{M}$ is surjective as $M$ is a generator in
${_{B}\mathcal{M}}$. Also, since $\mathrm{ev}_{M}$ is the same as the counit of
adjunction $(\mathbb{L}_B,\mathbb{R}_B)$,  $\mathrm{ev}
_{M}=\varepsilon_{B}^{B}$, the map $\mathrm{ev}_{M}$ is in the class $\eclass{M,B}$.
Hence the above sequence is  $\eclass{M,B}$-admissible and, for any $B$-bimodule
$N$, gives rise to
a long exact sequence, a part of which is:
\[
\mathrm{Ext}_{\eclass{M,B}}^{1}\!\left(  M\! \otimes_{A}\!{^{\ast}}\!M,N\right)\!
\rightarrow\mathrm{Ext}_{\eclass{M,B}}^{1}\!\left(  L,N\right)\!
\rightarrow\mathrm{Ext}_{\eclass{M,B}}^{2}\!\left(  B,N\right)\!
\rightarrow\mathrm{Ext}_{\eclass{M,B}}^{2}\!\left(  M\!\otimes_{A}\!{^{\ast}
}\!M,N\right)  .
\]
By Theorem~\ref{teo 4.2.33}, $M\otimes_{A}{^{\ast}}\!M={\mathbb{L}}_{B}\left(
{^{\ast}}M\right)  $ is $\eclass{M,B}$-projective so that $\mathrm{Ext}
_{\eclass{M,B}}^{1}\left(  L,N\right)
\simeq\mathrm{Ext}_{\eclass{M,B}}^{2}\left(  B,N\right)
 =\mathrm{H}_{M}^{2}(B,N)$. Hence the $\eclass{M,B}$-projectivity of
 $L=\ker (\mathrm{ev}_M)$ is equivalent to the property (c).
\end{proof}

Examples of formally smooth bimodules can be obtained from smooth extensions.
\begin{definition}\label{def.fse}
Let $A,$ $B$ be rings and let $i:A\rightarrow B$ be a ring homomorphism.
Consider the adjunction
\[
\begin{tabular}
[c]{lll}
$\mathbb{T}:{_{A}\mathcal{M}_{A}}\rightarrow{_{B}\mathcal{M}_{B},}$ &  &
${\mathbb{H}}:{_{B}\mathcal{M}_{B}}\rightarrow{_{A}\mathcal{M}_{A}}$\\
$\mathbb{T}\left(  X\right)  =B\otimes_{A}X\otimes_{A}B{,}$ &  & ${\mathbb{H}
}\left(  Y\right)  ={Y.}$
\end{tabular}
\ \ \ \ \ \ \
\]
Then $i$ is called a \emph{formally smooth extension} whenever $\mathrm{Ker}
\left(  m_{B}\right)  $ is $\mathcal{E}_{{\mathbb{H}}}$-projective. Here $m_B: B\otimes_A B\to B$ is the multiplication map and $\mathcal{E}_{{\mathbb{H}}}$
is a class of $\H$-relatively split morphisms as in Theorem~\ref{teo 4.2.33}.
\end{definition}
By \cite[Corollary 3.12]{AMS1}, ring extension $A\to B$ is formally
smooth  provided $B$ is formally smooth when
regarded as an algebra in the monoidal category $\left(  {_{A}\mathcal{M}
_{A},\otimes}_{A}{,A}\right)  .$

\begin{lemma}
\label{lem: ker mS} Let $A,$ $B$ be rings and let  $M$ be a $B$-$A$
bimodule  that is  finitely generated and
projective as a left $B$-module. Let $S={_{B}\mathrm{End}\left( M\right) }
=\mathbb{R}_{A}\left(  M\right)  \simeq{^{\ast}\!M\otimes
}_{B}M$ be the endomorphism ring. Write $m_S: S\otimes_A S\to S$
for the multiplication map and $\left(  L,i_{L}\right)$ for the kernel of
$\mathrm{ev}_M$. 
The following sequence
\[
\xymatrix{
0\ar[r] & {}^*\!M\ot_BL\ot_B M \ar[rr]^{\hspace{15pt}{}^*\!M\ot_Bi_L\ot_B M } &&
S\otimes_{A}S \ar[r]^{
\hspace{15pt}m_S} & S\ar[r]& 0.}
\]
is exact.
\end{lemma}

\begin{proof}
Start with the exact sequence
\[
\xymatrix{
0\ar[r] & L \ar[r]^{i_L\hspace{20pt}} & M\otimes_{A}{^{\ast}}\! M \ar[r]^{
\hspace{15pt}\mathrm{ev}_{M}} & B.} \eqno{(*)}
\]
Since $M$ is a finitely generated and projective left $B$-module,
 ${^{\ast}\!M}$ is a finitely generated and projective right
$B$-module. By tensoring ($*$) on the left with ${}^{\ast}\!M$
and on the right with $M$, we obtain the following exact sequence
\[
\xymatrix{
0\ar[r] & {}^*\!M\!\ot_B\!L\!\ot_B\! M \ar[rr]^{{}^*\!M\ot_Bi_L\ot_B M\hspace{15pt} } &&
{^{\ast}\!M \otimes}_{B} M \otimes
_{A} \!{^{\ast}}\!M {\otimes}_{B} M\;\; \ar[rr]^{\hspace{15pt}{^{\ast}\!M \otimes}_{B}\mathrm{ev}
_{M}{\otimes}_{B}M} &&\;{^{\ast}\! M\otimes}_{B}B{\otimes}_{B}M.}
\]
The isomorphisms ${^{\ast}M\otimes}_{B}B{\otimes}_{B}M\simeq{^{\ast
}M\otimes}_{B}M\simeq S,$  allow one to identify  the map
${^{\ast}M\otimes}_{B}\mathrm{ev}
_{M}{\otimes}_{B}M$  with  $m_{S}$. Being a multiplication of unital rings
the latter is
surjective.
\end{proof}

\begin{proposition}
\label{pro: FSBimod} Let $A,$ $B$ be rings and let  $M$ be a $B$-$A$
bimodule  that is  finitely generated and
projective as a left $B$-module. Let $S={_{B}\mathrm{End}}\left( M\right)
  =\mathbb{R}_{A}\left(  M\right)  \simeq{^{\ast}\!M\otimes
}_{B}M$ be the endomorphism ring. Write $i$ for the canonical ring map
$$
i: A\to S, \qquad a\mapsto [m\mapsto ma].
$$
\begin{zlist}
\item  If the bimodule ${}_BM_A$ is formally smooth, then
$i:A\rightarrow S$ is
a formally smooth extension.
\item If $M$ is a generator in ${_{B}\mathcal{M}}$ and
$i: A\rightarrow S$ is
a formally smooth extension, then ${}_BM_A$ is a formally smooth bimodule.
\end{zlist}
\end{proposition}

\begin{proof}
(1) In view of Theorem~\ref{teo 4.2.33}, to prove that $i:A\to S$
is a formally smooth extension, suffice it to prove that
$\mathrm{Ker}\left(  m_{S}\right)$ is a direct summand (in ${_{S}\mathcal{M}_{S}}$) of
$S\otimes_{A}X\otimes_{A}S$,
for a suitable object $X\in{_{A}\mathcal{M}_{A}.}$ Write $(L,i_L)$ for the kernel
of $\mathrm{ev}_M$. Since $M$ is formally smooth, $L$ is $\eclass{M,B}$-projective.
By Theorem~\ref{teo 4.2.33}, this
means that the counit of the adjunction $\left(  {\mathbb{L}}_{B},{\mathbb{R}
}_{B}\right)  $ evaluated at $L$
\[
\varepsilon_{L}^B:{\mathbb{L}}_{B}{\mathbb{R}}_{B}\left(  L\right)  \simeq
M\otimes_{A}{^{\ast}}M\otimes_{B}L\rightarrow L
\]
has a section $\sigma:L\rightarrow M\otimes_{A}{^{\ast}}M\otimes_{B}L$ in
${_{B}\mathcal{M}_{B}.}$

Since $M$ is a finitely generated and projective left $B$-module, the functor
$\mathbb{R}_B$ can be naturally identified with the tensor functor
${}^*\! M\otimes_B -$. Furthermore, $\mathrm{ev}_M$ is in the class $\eclass{M,B}$,
hence ${^{\ast}}\!M\otimes_{B}\mathrm{ev}_{M}\simeq{\mathbb{R}}_{B}\left(
\mathrm{ev}_{M}\right)$
splits in ${_{A}\mathcal{M}_{B}}$. Thus applying $\mathbb{R}_B$ to the defining sequence
of $(L,i_L)$ we obtain the split exact sequence of $A$-$B$ bimodules
\[
\xymatrix{0\ar[r] & {^{\ast}}\! M\otimes_{B}L\ar[rr]^-{{^{\ast}}M\otimes_{B}i_{L}
} && {^{\ast}}\!M\otimes_{B}M\otimes_{A}{^{\ast}}\!M \ar[rr]^-{{^{\ast
}}\!M\otimes_{B}\mathrm{ev}_{M}} & & {^{\ast}}\! M \ar[r] & 0 }
\]
In particular ${}^*M\otimes_Bi_L$ is a section in ${}_A\mathcal{M}_B$, and, consequently
$M\otimes_A{}^*M\otimes_Bi_L$ is a section in ${}_B\mathcal{M}_B$. Therefore, the map
\[
\alpha:=\left(  \xymatrix{L \ar[r]^-{\sigma} & M\otimes_{A}{^{\ast}
}\!M\otimes_{B}L \ar[rr]^-{M\otimes_{A}{^{\ast}}M\otimes_{B}i_{L}}
&& M\otimes_{A}{^{\ast}}\!M\otimes_{B}M\otimes_{A}{^{\ast}}\!M}\right)
\]
splits in ${_{B}\mathcal{M}_{B}.}$ Consequently,
\[\xymatrix{
{^{\ast}\!M\otimes}_{B}L{\otimes}
_{B}M\ar[rr]^-{{^{\ast}M\otimes}_{B}\alpha{\otimes}_{B}M} && {^{\ast}\!M\otimes}_{B}M\otimes_{A}{^{\ast}}\!M\otimes_{B}M\otimes_{A}{^{\ast}
}\!M{\otimes}_{B}M\simeq S\otimes_{A}S\otimes_{A}S}
\]
splits in ${_{S}\mathcal{M}_{S}}$. In view of Lemma~\ref{lem: ker mS},
$\mathrm{Ker}\left(  m_{S}\right)  \simeq{^{\ast}M\otimes}_{B}L{\otimes}
_{B}M$, and hence ${^{\ast}M\otimes}_{B}\alpha{\otimes}_{B}M$
 is the required $S$-bimodule section
of a map $S\otimes_A S\otimes_A S\to \ker (m_S)$.

(2) By \cite[Theorem~3.8 \& Theorem~4.42]{AMS1}, if $i:A\to S$ is formally
smooth,
then  $\mathrm{Hdim}(S|A)\leq 1$. Since $M$ is a generator left $B$-module,
Theorem~\ref{teo: Hoch S} implies that $\mathrm{Hdim}_M(B) \leq 1$.
Proposition~\ref{prop.smooth} then implies that $M$ is a formally smooth bimodule.
\end{proof}

\begin{proposition}\label{prop.smooth.alg}
Let $B$ be an algebra over a commutative ring $k$. Consider the functor
\[
\mathbb{F}:{\mathcal{M}_{k}}\rightarrow{_{k}\mathcal{M}_{k}
}:\left(  V,\mu^{r}\right)  \longmapsto\left(  V,\mu^{l},\mu^{r}\right)  ,
\]
where the left $k$-multiplication is defined by $\mu^{l}  \left(
\lambda\otimes _k v \right) : =\mu^{r} \left( v\otimes_k\lambda
\right)$, for all $\lambda \in k$ and $v\in V$. Furthermore,
consider the adjunction
\[
\begin{tabular}
[c]{lll}
$\mathbb{T}^{\prime}:{\mathcal{M}_{k}}\rightarrow{_{B}\mathcal{M}_{B},}$ &
& ${\mathbb{H}}^{\prime}:{_{B}\mathcal{M}_{B}}\rightarrow{\mathcal{M}_{k}
}$\\
$\mathbb{T}^{\prime}\left(  X\right)  =B\otimes_{k}\mathbb{F}\left(
X\right)  \otimes_{k}B{,}$ &  & ${\mathbb{H}}^{\prime}\left(  Y\right)
={Y.}$
\end{tabular}
\ \ \ \ \ \
\]
The following assertions are equivalent.
\begin{blist}
\item The bimodule $_{B}M_{k}={_{B}B_{k}}$ is formally smooth.
\item $\mathrm{Ker}\left(  m_{B}\right)  $ is $\mathcal{E}_{{\mathbb{H}
}^{\prime}}$-projective.
\item The extension $k\rightarrow{B}$ is formally smooth.
\end{blist}
\end{proposition}
\begin{proof}
Clearly {the} $B$-module $B$ can be identified both with its
dual and endomorphism ring. With this identification the
evaluation map $\mathrm{ev}_B = m_B$. Hence the equivalence (a) $
\Leftrightarrow $ (c) follows by Proposition~\ref{pro: FSBimod}.
The implication (b) $\Rightarrow$ (c) is an immediate consequence
of the observation that any $B$-bimodule map that splits as a
$k$-bimodule map, splits also as a right
$k$-module map (i.e.\ $\mathcal{E}
_{{\mathbb{H}}}\subseteq\mathcal{E}_{{\mathbb{H}}^{\prime}}$, where $\mathbb{H}$ is
a functor in Definition~\ref{def.fse}).

(a) $\Rightarrow$ (b) We need to show that $L:=\mathrm{Ker}\left(
\mathrm{ev}_{{B}}\right)  $ is
$\mathcal{E}_{{\mathbb{H}}^{\prime}}$-projective. The counit of
the adjunction $(\mathbb{T}^{\prime}, {\mathbb{H}}^{\prime})$ is
given by the two-sided multiplication
\[
\varepsilon'_{N}:B\otimes_{k}\mathbb{F}{\mathbb{H}}^{\prime}\left(
N\right)  \otimes_{k}B\rightarrow N,\text{ for every }N\in{_{B}
\mathcal{M}_{B}}.
\]
Note that $\varepsilon'_{N}\in\mathcal{E}_{{\mathbb{H}}^{\prime}}$ as it is the
counit. Consider the adjunction
\[
\begin{tabular}
[c]{lll}
$\mathbb{T}^{\prime\prime}:{\mathcal{M}_{k}}\rightarrow{\mathcal{M}_{B},}$
&  & ${\mathbb{H}}^{\prime\prime}:{\mathcal{M}_{B}}\rightarrow{\mathcal{M}
_{k}}$\\
$\mathbb{T}^{\prime\prime}\left(  X\right)  =X\otimes_{k}B{,}$ &  &
${\mathbb{H}}^{\prime\prime}\left(  Y\right)  ={Y.}$
\end{tabular}
\ \ \ \ \ \
\]
By the standard argument (cf.\ e.g.\ \cite[Proposition~2.5]{CunQui:alg}),
$L\simeq\frac{B}{k}\otimes_{k
}B=\mathbb{T}^{\prime\prime}\left(  B/k\right)$. The latter is $\mathcal{E}
_{\mathbb{H}^{\prime\prime}}$-projective by Theorem~\ref{teo
4.2.33}. Since
$\varepsilon'_{L}\in\mathcal{E}_{{\mathbb{H}}^{\prime}}\subseteq$
$\mathcal{E}_{\mathbb{H}^{\prime\prime}}$ we conclude that
$\varepsilon'_{L}$ splits in ${\mathcal{M}_{B}}$, that is
$\varepsilon^{\prime}_{L}\in\eclass{M,B}$ (note that
$L=\mathbb{F}{\mathbb{H}}^{\prime}\left(  L\right)  $ as it is a
subbimodule of $B\otimes_{k}B$). By hypothesis $L$ is
$\eclass{M,B}$-projective so that $\varepsilon'_{L}$ splits in
${_{B}\mathcal{M}_{B}.}$ By Theorem~\ref{teo 4.2.33}~(c)$
\Rightarrow$(a), we thus conclude that $L$ is $\mathcal{E}_{{\mathbb{H}
}^{\prime}}$-projective.
\end{proof}

In view of Proposition~\ref{prop.smooth.alg} a formally smooth
algebra $B$ over a field $k$ is a formally smooth
$(B,k)$-bimodule. In this way one can construct examples of
formally smooth bimodules which are not separable. For example,
the tensor algebra $T_k(V)$ of a vector space $V$ is formally
smooth but not separable in view of Proposition \ref{pro: sep S}.
In fact it is well known that any separable extension of a field
$k$ is finite dimensional over $k$ (cf. \cite[{Proposition
1.1}]{Villamayor}).

Let $M$ be a left $B$-module and write $S$ for its endomorphism ring.
Recall that a left $B$-module $N$ is said to be
{\em
$M$-static} provided  the
evaluation
\[
\mathrm{ev}_{M,N}:M\otimes_{S}{_{B}\mathrm{Hom}}\left( M ,N\right)
  \rightarrow N, \qquad m\otimes_A f\mapsto (m)\! f
\]
is an isomorphism  (see e.g.\ \cite[2.3]{Wisbauer}).  Recall further that
the image of the evaluation map $\mathrm{ev}_{M,N}$,
is called the {\em trace of $M$ in $N$} and is denoted by $\mathrm{Tr}_M(N)$.
Finally,
denote by $\sigma[M]$ the full subcategory of ${}_B\mathcal{M}$, whose objects
are all modules subgenerated by $M$ (cf.\ \cite[Section~15]{Wisbauer-book}).

\begin{proposition}
\label{pro: M-static} Let $B$ be a ring,  ${M}$ be a left $B$-module and set
$S={_{B}\mathrm{End}}\left(M\right)$, so that $M$ is a $B$-$S$ bimodule.
The  following assertions are equivalent:
\begin{blist}
\item The evaluation map
$\mathrm{ev}_{M}:M\otimes_{S} {}^*\!M  \rightarrow B$ is injective.
\item The $B$-module $\mathrm{Tr}_{M}\left(  B\right)  $ is $M$-static.
\end{blist}

In particular these conditions hold whenever $M$ is a generator in $\sigma\left[  M\right]  $.
\end{proposition}

\begin{proof}
Since
\[
{_{B}\mathrm{Hom}}\left(  M,\mathrm{Tr}_{M}\left(  B\right)  \right)
={_{B}\mathrm{Hom}}\left( M,B\right)  ,
\]
the equivalence follows by observing that
$\mathrm{ev}_{M,\mathrm{Tr}_{M}\left(  B\right)}$ is exactly $\mathrm{ev}_{M}$
corestricted to its image.  The last assertion follows by
\cite[Lemma~1.3]{Zim:end}.
\end{proof}

Combining Proposition~\ref{pro: M-static} with Proposition~\ref{prop.1ex}
we immediately obtain
\begin{corollary} \label{cor.sigma}
If a left $B$-module $M$  with endomorphism ring $S$ generates $\sigma[M]$, then
$M$ is a formally smooth
$B$-$S$ bimodule.
\end{corollary}

\begin{corollary}
\label{coro: M-static} Let $B$ be a ring,  $M$ be a left $B$-module, and let $S={_{B}\mathrm{End}}\left(  M\right)  $.
The following assertions are equivalent:

\begin{blist}
\item ${_{B}}M_{S}$ is a separable bimodule.

\item $M$ is a generator in ${_{B}\mathcal{M}}$.

\item The evaluation map $\mathrm{ev}_{M}:M\otimes_{S} {}^*\!M  \rightarrow B$
is an isomorphism.
\end{blist}
\end{corollary}

\begin{proof}
The implication (a) $\Rightarrow$ (b) follows by Proposition~ \ref{pro: SepBimod}.
If  $M$ is a generator of ${}_B\mathcal{M}$, then it is also a generator of
$\sigma\left[M\right]  $. By Proposition~\ref{pro: M-static},
$B=\mathrm{Tr}_{M}\left(  B\right)  $ is $M$-static. Hence $\mathrm{ev}_{M}$
is an isomorphism. This proves that (b) implies (c). The implication
(c) $\Rightarrow$ (a)
is obvious.
\end{proof}

The following proposition explains how two formally smooth
bimodules can be combined to give a formally smooth bimodule,
and thus can be seen as module version of
\cite[Proposition~5.3]{CunQui:alg}.

\begin{proposition}\label{prop.constr}
Let $A,$ $B$ and $T$ be rings. Let  $Y$ be a $T$-$A$ bimodule and let
$X$ be a $B$-$T$ bimodule
such that
 the evaluation map $\mathrm{ev}_{X}: X\ot_T{}^*\! X\to B$ is injective and that
$X$ is  flat as a right $T$-module.
Assume that one of the following conditions (1) or (2) is satisfied:
\begin{zlist}
\item$Y$ is a separable $T$-$A$-bimodule.
\item
\begin{rlist}
\item  ${^{\ast}\!X}  $ is flat as a right $T$-module,
\item $Y$ is  finitely generated and projective as a left $T$-module, and
\item $Y$ is a formally smooth $T$-$A$ bimodule.
\end{rlist}
\end{zlist}
Then
\[
{_{B}M_{A}}={_{B}}X\otimes_{T}Y_{A}.
\]
is a formally smooth bimodule.

In particular, if a left $B$-module $X$ is a generator of
$\sigma\left[ X\right]  $, $T={_{B}\mathrm{End}}\left(
X\right)  $, and either (1) or (2) hold, then $M=X\otimes_{T}Y$
is a formally smooth $B$-$A$
bimodule.
\end{proposition}

\begin{proof} Associate with $Y$ the tensor-hom adjunction,
\[
\begin{tabular}
[c]{lll}
$\mathbb{T}:{_{A}\mathcal{M}_{B}}\rightarrow{_{T}\mathcal{M}_{B},}$ &
& ${\mathbb{H}}:{_{T}\mathcal{M}_{B}}\rightarrow{_{A}\mathcal{M}_{B}}
$\\
$\mathbb{T}\left(  U\right)  =Y\otimes_{A}U{,}$ &  & $\mathbb{H}
\left(  W\right)  ={_{T}\mathrm{Hom}}\left(  Y,
W\right)  {,}$
\end{tabular}
\ \ \
\]
and denote its counit (the evaluation map) by $\varepsilon$.
 Use the natural isomorphism
\[
\Phi: {^{\ast}M}
={_{B}\mathrm{Hom}}\left(  X\otimes_{T}Y ,B\right)
\to {_{T}\mathrm{Hom}}\left( Y
,{}^{\ast}\!X\right), \qquad f\mapsto [y\mapsto f(-\ot_Ty)],
\]
to write the evaluation map $\mathrm{ev}_M: M\ot_A {}^*\!M\to
{B}$  as
\[
\mathrm{ev}_{M}=\mathrm{ev}_{X}\circ\left(  X\otimes_{T}\varepsilon
 _{{^{\ast}\!X}}\right)  \circ\left(  M\otimes_{A}
\Phi\right)  .
\]
Since, by assumption, $X_T$ is flat and $\mathrm{ev}_X$ is injective, and since
 $\Phi$ is an isomorphism, there is an isomorphism of $B$-bimodules
\[
\mathrm{Ker}\left(  \mathrm{ev}_{M}\right)  \simeq X\otimes_{T}\mathrm{Ker}\left(  \varepsilon
 _{{^{\ast}\!X}}\right)  . \eqno{(*)}
\]

Assume that condition (1) is satisfied, i.e.\ that $Y$ is a separable bimodule.
By Proposition~\ref{pro: SepBimod},
${\mathbb{H}}$ is a separable functor, hence,  by Corollary~\ref{coro Rafael},
any object in ${_{T}\mathcal{M}_{B}}$ is $\mathcal{E}_{\mathbb{H}}
$-projective. In particular
$\ker\left( \varepsilon_{{^{\ast}\!X}}\right)$ is $\mathcal{E}_{\mathbb{H}}$-projective.
 By Theorem~\ref{teo 4.2.33}, $\ker\left( \varepsilon_{{^{\ast}\!X}}\right)$
  is a
direct summand in ${_{T}\mathcal{M}_{B}}$ of $\mathbb{T}\left(
U\right)  $ for some $U\in{_{A}\mathcal{M}_{B},}$ hence
$X\otimes
_{T} \ker\left( \varepsilon_{{^{\ast}\!X}}\right)$ is a direct summand of
\[
X\otimes_{T}\mathbb{T}\left(  U\right)  =X\otimes_{T}Y\otimes
_{A}U=M\otimes_{A}U=\mathbb{L}_{B}\left(  U\right)
\]
in ${_{B}\mathcal{M}_{B}.}$ Theorem~\ref{teo 4.2.33} implies that
$\mathrm{Ker}\left(  \mathrm{ev}_{M}\right)  $ is
$\eclass{M,B}$-projective so that $M$ is a formally smooth bimodule.

Assume that conditions (2) hold.
Since $Y$ is a finitely generated and projective left $T$-module, the functor
$\mathbb{H}$ is naturally isomorphic to the tensor functor
${}^*Y\otimes_T -$, and the counit $\varepsilon$ evaluated at
$W$ can be identified with $\mathrm{ev}_Y\ot_T W$. In particular,
$\ker \left(\varepsilon_{{}^*\!X}\right) \simeq \ker\left(\mathrm{ev}_Y\ot_T {}^*\!X\right)$. Since ${}^*\!X$ is a flat left $T$-module, the isomorphism ($*$) yields
\[
\mathrm{Ker}\left(  \mathrm{ev}_{M}\right)  \simeq
X\otimes_{T}\mathrm{Ker\,}\left(
\mathrm{ev}_{Y}\right)  \otimes_{T}{^{\ast}X}.
\]
Since ${_{T}}Y_{A}$ is a formally smooth bimodule, $\mathrm{Ker}
\left(  \mathrm{ev}_{Y}\right)$ is $\mathcal{E}_{\mathbb{H}}
$-projective, which, as in the case (1), implies that $M$ is a formally smooth
$B$-$A$ bimodule.

To prove the final statement observe that if $_{B}X$ is a
generator of $\sigma\left[ X\right]  $ and $T={_{B}\mathrm{End}}\left(
X\right)  ,$ then by Proposition~\ref{pro: M-static}, $\mathrm{ev}_{X}$
is injective. Furthermore by \cite[Section~15.9]{Wisbauer-book}, $X_{T}$ is
a flat, hence the main assumptions about $X$ are satisfied.
\end{proof}

In \cite[Section~2]{CaeZhu:sep} several ways of constructing separable bimodules
are described. Combined with Proposition~\ref{prop.constr} these can
 provide a source of examples of smooth bimodules.

\end{document}